\documentclass[12pt]{article}
\usepackage{amsmath, amssymb, amsthm, color, comment}
\usepackage{tikz, here}
\usetikzlibrary{cd}
\usepackage{hyperref}
\usepackage{cleveref}

        \usepackage{setspace}

\usepackage[margin=1in]{geometry}
\hypersetup{
    colorlinks=true,
    linkcolor=blue,
    filecolor=magenta,      
    urlcolor=cyan,
}
\newtheorem{theorem}{Theorem}[section]
\newtheorem{corollary}[theorem]{Corollary}
\newtheorem{lemma}[theorem]{Lemma}

\newtheorem{proposition}[theorem]{Proposition}
\newtheorem{claim}{Claim}

\usepackage{etoolbox}
\AtEndEnvironment{theorem}{\setcounter{claim}{0}}
\AtEndEnvironment{lemma}{\setcounter{claim}{0}}
\AtEndEnvironment{corollary}{\setcounter{claim}{0}}
\AtEndEnvironment{proposition}{\setcounter{claim}{0}}

\theoremstyle{remark}

\theoremstyle{definition}
\newtheorem{definition}[theorem]{Definition}

\renewenvironment{proof}[1][Proof.]{\begin{trivlist}
\item[\hskip \labelsep {\textit{#1}}]}{\end{trivlist}}

\newcommand{\ult}{\text{Ult}}
\newcommand{\crit}{\text{crit}}

\newcommand{\ran}{\text{ran}}

\newcommand{\hull}{\text{Hull}}

\newcommand{\restrict}{\!\upharpoonright\!}
\newcommand{\is}{\trianglelefteq}

\newcommand{\tree}[1]{\mathcal{#1}}

\makeatletter
\tikzset{
  edge node/.code={%
      \expandafter\def\expandafter\tikz@tonodes\expandafter{\tikz@tonodes #1}}}
\makeatother
\tikzset{
  is/.style={
    draw=none,
    edge node={node [sloped, allow upside down, auto=false]{$\is$}}},
  Is/.style={
    draw=none,
    every to/.append style={
      edge node={node [sloped, allow upside down, auto=false]{$\is$}}}
  }
}
\tikzset{
  eq/.style={
    draw=none,
    edge node={node [sloped, allow upside down, auto=false]{$=$}}},
  Eq/.style={
    draw=none,
    every to/.append style={
      edge node={node [sloped, allow upside down, auto=false]{$=$}}}
  }
}

\title{The uniqueness of the core model}
\author{Benjamin Siskind\footnote{This research was funded in whole or in part by the Austrian Science Fund (FWF) [10.55776/Y1498]. For open access purposes, the author has applied a CC BY public copyright license to any author accepted manuscript version arising from this submission.}\\
TU Wien\\
Vienna, Austria}

\begin{document}
\maketitle
\begin{abstract}
    The Jensen-Steel core model is a canonical inner model which plays a fundamental role in the meta-mathematics of set theory. Its definition depends on exactly which hierarchy of fine-structural models of set theory, premice, one uses. Each such hierarchy involves somewhat arbitrary decisions and working with different hierarchies ostensibly leads to different versions of the core model. We show that in some contexts, abstract properties of the core model uniquely determine it; that is, there is at most one inner model with these properties.
\end{abstract}

\section{Introduction}
A major achievement of inner model theory has been the identification of the core model $K$ by Ronald Jensen and John Steel in \cite{kwm},  under the hypothesis that there is no inner model with a Woodin cardinal. $K$ plays a central role in the meta-mathematics of set theory: it is an essential tool in establishing strong consistency strength lower bounds for natural theories, for example in Steel's result that PFA implies AD${}^{L(\mathbb{R})}$ \cite{pfa}. In this result and others like it, one uses that when there is no inner model with a Woodin cardinal, $K$ is a canonical inner model which is close to $V$ (a generalization of Jensen's seminal result that $L$ is close to $V$ under the more restrictive hypothesis that $0^\#$ doesn't exist).

Jensen and Steel identify $K$ as an inner model whose levels are certain \textit{premice}, fine-structural models of set theory which have a complicated definition. In particular, Jensen and Steel use what are known as ms-indexed (pure-extender) premice to build $K$. Other varieties of premice have been studied, for example Jensen-indexed premice, or the recent pfs-premice and least branch strategy mice of Steel's \cite{steel book}. Using these different varieties of premice could give rise to ostensibly different versions of $K$, however it is expected that all these versions are actually the same. One reason for this expectation is that it should be possible to translate premice of one variety into premice of another variety. This has been realized in some cases; for example, Fuchs \cite{fuchs1} and \cite{fuchs2} showed that one can translate ms-indexed premice into a modified Jensen-indexed premice, and vice-versa. These translation methods are carefully tailored to the varieties one is translating between, so such methods don't seem like they can yield the kind of general result one would really like to show: \textit{any} successful notion of premouse must give rise to the same core model. 

In this paper we take a new approach to establishing sufficiently general results along these lines. We show that in some contexts, abstract properties of the core model uniquely determine it; that is, there is at most one inner model with these properties. So any notion of premouse for which the associated core model enjoys the abstract properties of the known core model must actually give rise to the same model. We should mention that this kind of result is not without precedent. For example, Steel showed that $\text{HOD}^{L(\mathbb{R})}|\Theta^{L(\mathbb{R})}$ is the universe of an ms-indexed premouse, so that the universe of this premouse is identifiable without reference to the particular fine structure at all \cite{hod_as}. Moreover, some ideas at the heart of our proofs are essentially present in this kind of $\text{HOD}$ computation.

This paper is organized as follows. In \S3, we characterize the core model under the additional hypothesis that there is a proper class of measurable cardinals. In \S4, we'll provide a different characterization of the core model which does not need measurable cardinals, but needs the assumption that $0^\text{\textparagraph}$ does not exist. 
We start with some preliminary definitions and observations and some facts about $K$.

The work in this paper appeared in the author's PhD thesis \cite{thesis}. The author would like to thank John Steel for many useful discussions about the ideas in this paper.

\section{Preliminaries}
We will consider transitive models of ZFC${}^-$, that is ZFC, stated with the Replacement Schema and the Well-Ordering Theorem,\footnote{The Well-Ordering Theorem is the statement that every set admits a well-ordering.} but without the Power Set Axiom. We introduce the following bits of notation.
\begin{definition}
For $M$ a transitive model of ZFC${}^-$, $o(M)=\text{Ord}\cap M$. For $M, N$ transitive models of ZFC${}^-$, and $\pi:M\to N$ an elementary embedding, we let $\pi(o(M))=o(N)$.
\end{definition}
We also introduce the following nonstandard notation, for convenience.
\begin{definition}
For $\mu$ a limit cardinal we let \[H_\mu=\bigcup\{H_\kappa\mid \kappa<\mu \text{ a regular cardinal}\}.\]
For $M$ a transitive model of ZFC${}^-$, we also set $H_{o(M)}^M=M$.
\end{definition}
This notation is useful because this hierarchy comes up naturally in Inner Model Theory: for $M$ a premouse and $\mu$ a limit cardinal of $M$, $H_\mu$ is the universe of $M|\mu$.


We now review some well-known facts about $K$. First, we will state a folklore theorem about the absoluteness of iterability when there is no inner model with a Woodin cardinal.

\begin{theorem} \label{absoluteness}
Assume there is no inner model with a Woodin cardinal. Let $W$ be an inner model of ZFC. Then $K^W$ is iterable.
\end{theorem}

This is a corollary of the following result, due to Steel (this is easy to obtain from iterability absoluteness results in \cite{hod_as} and standard facts about the existence of $Q$-structures for normal iteration trees on 1-small premice).

\begin{theorem}\label{absoluteness2}
Let $W$ be an inner model of ZFC, $\kappa$ be an uncountable cardinal of $W$, and $P$ a 1-small premouse with $P\in H_\kappa^W$. Then $P$ is iterable iff $H_\kappa^W\models ``P$ is iterable".
\end{theorem}

\begin{proof}[Proof of Theorem \ref{absoluteness}.]
Since there is no inner model with a Woodin cardinal, $K^W$ is defined and also has no Woodin cardinals. It follows that for any successor $K^W$-cardinal, $\eta$, $K^W|\eta$ is 1-small. But if there is a bad normal iteration tree $\tree{T}$ on $K^W$ in $V$, then there is a bad tree on $K^W|\eta$ for some such $\eta$. But $K^W|\eta$ is iterable in $W$ and so iterable in $V$, by Theorem \ref{absoluteness2}, a contradiction. \qed 
 
\end{proof}

Next, we'll review Steel's inductive definition of $K$ from \cite{cmip}. We could likely also use Schindler's result from \cite{schindler} that, above $\omega_2$, levels of $K$ are just obtained by stacking, though this has not been checked in context without the measurable cardinal.

\begin{definition}
For $\alpha$ a $K$-cardinal, a countably iterable, 1-small premouse $N$ is \textit{$\alpha$-strong} iff $K|\alpha\is N$ and for all premice $M$ such that $M$ is $\beta$-strong for all $K$-cardinals $\beta<\alpha$, the phalanx $(N, M, \alpha)$ is iterable.
\end{definition}

Arguments from \cite{cmip} give
\begin{theorem} \label{definability}
Assume there is no inner model with a Woodin cardinal. Let $\alpha$ be a cardinal of $K$. Then 
\begin{enumerate}
    \item $N$ is $\alpha$-strong iff $K|\alpha\is N$ and for all premice $M$ of size $\leq |N|$ such that $M$ is $\beta$-strong for all $K$-cardinals $\beta<\alpha$, the phalanx $(M, N, \alpha)$ is $\omega_1$-iterable;
    \item $K|\alpha^{+,K}=\bigcup\{N|\alpha^{+,N}\mid N\text{ is $\alpha$-strong and $|N|=|\alpha|$}\}$.
\end{enumerate}
\end{theorem}

This immediately gives the following result.

\begin{theorem}\label{definability of K 1}
Assume there is no inner model with a Woodin cardinal. Let $\mu$ be a strong limit cardinal. Then $K|\mu$, i.e. $H^K_\mu$ together with the extender sequence of $K|\mu$, is definable without parameters over $H_{\mu}$, uniformly in $\mu$.
\end{theorem}

\begin{proof}
Fix $\mu$. The idea here is to define when some premouse is an initial segment of $K|\mu$ by asserting there are sufficiently long sequences of sets $S_\alpha$, ordinals $\kappa_\alpha$, and premice $P_\alpha$ such that the $\kappa_\alpha$ are the $K$ cardinals, $S_\alpha$ is the set of $<\kappa_\alpha$-strong premice of some \textit{fixed} cardinality, and $P_\alpha=K|\kappa_\alpha$, using the inductive definition of $K$ (i.e. the previous theorem). We use that $\mu$ is a strong limit to guarantee that the set of all  $<\kappa_\alpha$-strong premice of our fixed cardinality $<\mu$ is a member of $H_\mu$, since premice of size $\theta$ are essentially subsets of $\theta$.

This is routine, but we include it here for completeness. A premouse $Q$ of size $<\mu$ is a proper initial segment of $K|\mu$ iff there sequences $\langle S_\alpha\mid \alpha\leq\zeta\rangle$,  $\langle P_\alpha\mid \alpha\leq\zeta\rangle$, $\langle \kappa_\alpha\mid \alpha\leq\zeta\rangle$, for some $\zeta<\mu$, such that 
\begin{enumerate}
    \item \begin{itemize}
        \item\ $S_0$ is the set of countably iterable, 1-small premice of size $\leq |Q|$, 
         \item $P_0=\langle V_\omega,\emptyset\rangle$,and 
         \item $\kappa_0=\omega$;
         \end{itemize}
    \item for $\alpha+1\leq\zeta$,
    \begin{itemize}
        \item $S_{\alpha+1}$ is the set of all $N\in S_\alpha$ such that $P_\alpha\is N$ and for all $M\in S_\alpha$, $(M, N, \alpha)$ is $\omega_1$-iterable,
        \item $P_{\alpha+1}=\bigcup\{N|\kappa_\alpha^{+,N}\mid N\in S_{\alpha+1}\}$, and
        \item $\kappa_{\alpha+1}=o(P_{\alpha+1})$;
    \end{itemize}
    \item for $\lambda\leq\zeta$ a limit ordinal,
    \begin{itemize}
        \item $S_\lambda=\bigcup \{S_\alpha\mid \alpha<\lambda\}$,
        \item $P_\lambda=\bigcup \{P_\alpha\mid \alpha<\lambda\}$, and
        \item $\kappa_\lambda=\sup\{\kappa_\alpha\mid \alpha<\lambda\}$; and
    \end{itemize}
\item $Q\is P_\zeta$.
\end{enumerate}
By Theorem \ref{definability} and our above comments, this gives a definition for $K|\mu$ over $H_\mu$ and is clearly uniform in $\mu$. 
\qed
\end{proof}

Fix $\varphi_K(v)$ the formula in the language of set theory which defines $K$ as in the previous theorem. We also let ``V=K" be the sentence $\forall v\, \varphi_K(v)$.

Theorem \ref{definability} also gives the following.

\begin{theorem} \label{k of k}
Assume there is no inner model with a Woodin cardinal. Then $K\models ``V=K"$. In particular, $K\models ``V=HOD"$.
\end{theorem}
\begin{proof}[Proof sketch.]
The point is that, by induction, we'll have that for all $K$-cardinals $\beta<\alpha$, $K|\alpha$ is $\beta$-strong inside of $K$. For $1$-small $N$ which is $\alpha$-strong in $K$, the iterability of $(K|\alpha, N, \alpha)$ inside $K$ implies that this phalanx is actually iterable in $V$, which suffices for showing that $N$ is actually $\alpha$-strong.

This implies $K\models``V=HOD"$ because $K$ has a global well-order definable over $V$, by Theorem \ref{definability of K 1}.\qed
\end{proof}

\noindent Theorem \ref{k of k} can be proved in other ways; for example, using Theorem \ref{k is local k}, below.

Unsurprisingly, one of the key properties we will use in identifying the core model is covering. Informally, the covering properties of an inner model $W$ are thought of asserting that $W$ is ``close" to $V$. The specific covering property will make use of in most of our results is captured in the following definition.
\begin{definition} For inner models $M\subseteq  N$ of ZFC, we say that $M$ is \textit{close} to $N$ if for all measurable or singular strong limit cardinals $\mu$ of $N$, 
\begin{enumerate}\item$\mu$ is measurable or singular in $M$ and \item $\mu^{+,M}=\mu^{+,N}$. 
\end{enumerate}
We'll say that $M$ is \textit{close} if $M$ is close to $V$.
\end{definition}

The Jensen-Steel core model $K$ is close, provably in ZFC+``there is no inner model with a Woodin cardinal". This follows by combining the covering theorems of from Jensen-Steel \cite{kwm} and Mitchell-Schimmerling \cite{mitchell-schimmerling}. That is, we have the following.
\begin{theorem}\label{k is close}
Assume there is no inner model with a Woodin cardinal. Then $K$ is a close inner model.
\end{theorem}

As far as we can tell, ordinary weak covering, i.e. that $\text{cof}(\lambda)\geq |\lambda|$ whenever $\lambda\geq\omega_2$ is a successor cardinal of $M$, may not be transitive whereas the property just introduced is transitive, that is we have the following.\footnote{There is a consequence of weak covering which is transitive: if $\mu\geq \omega_2$ is a regular cardinal, then $cof(\mu^{+,M})\geq \mu$. This actually works fine for our purposes below $0^\text{\textparagraph}$ but does not seem to work below a Woodin cardinal, in general. 
}

\begin{proposition}\label{transitivity}
Suppose $M\subseteq N\subseteq P$ are inner models of ZFC, $M$ is close to $N$, and $N$ is close to $P$. Then $M$ is close to $P$.
\end{proposition}

\section{With measurable cardinals}
In this section we prove a uniqueness theorem about the core model under the hypothesis that there is no inner model with a Woodin cardinal, assuming that there is a proper class of measurable cardinals. We'll define what it means for an inner model to ``resemble the core model", show there is at most one such inner model, and then prove that $K$ is that model. 

As mentioned in the introduction, the problem motivating such a result is whether different notions of premice give rise to the same core model. A natural question is whether the core model associated to Jensen-indexed premice also resembles the core mode, in our sense. Unfortunately, we do not know whether this is the case at present. The problem is that the theory developed in \cite{cmip} or \cite{kwm} has not been fully worked out for Jensen-indexed premice. However, Jensen's in-progress manuscript \cite{jensen manu} will include the development of such a theory. We believe that the resulting core model will indeed satsify our definition and so our theorem would show that this Jensen-indexed core model has the same universe as the ms-indexed core model, $K$. 

Under the hypothesis that there is no inner model with a Woodin cardinal but there is a proper class of measurable cardinals, the ms-indexed $K$ is just $\bigcup\{K_\mu\mid \mu\text{ measurable}\}$, where $K_\mu$ is the core model from Steel's \cite{cmip} at the measurable cardinal $\mu$. In this context, we can identify $K$ in a particularly simple to state manner. For $\mu$ a measurable cardinal, we'll identify $H^K_{\mu^+}$ as the unique $H_{\mu^+}^M$ for $M$ which ``resembles the core model at $(\mu,\lambda)$", for $\lambda$ any inaccessible cardinal above $\mu$.

\begin{definition}
Let $\mu<\lambda$ with $
\mu$ measurable and $\lambda$ inaccessible. A transitive model $P$ is \textit{$\mu$-full at $\lambda$} iff $P=V_\lambda^W$ for an inner model $W$ of ZFC such that $\mu$ is measurable in $W$ and $\mu^{+,W}=\mu^+$.
\end{definition}

It is not immediately obvious that being $\mu$-full at $\lambda$ is expressible in the language of set theory, since we quantified over the proper class $W$ in the above definition. We leave it to the reader to check the following easy proposition.

\begin{proposition}
$P$ is $\mu$-full at $\lambda$ iff $P$ is a transitive model of ZFC, $o(P)=\lambda$, $P\models \mu$ is measurable, $\mu^{+,P}=\mu^+$, and there is a well-order $\leq$ of $P$ such that every bounded subset of $\lambda$ constructible from $P$ and $\leq$ is in $P$ (i.e. $P=V^{L(P,\leq)}_\lambda$).
\end{proposition}
 We let $\text{Full}_{\mu,\lambda}$ be the set of all $P$ which are $\mu$-full at $\lambda$.
We can now state the main definition of the section.
\begin{definition}
Let $\mu$ be a measurable cardinal and $\lambda>\mu$ inaccessible. A transitive model $M$ \textit{resembles the core model at $(\mu,\lambda)$} if there is a function from $\text{Full}_{\mu,\lambda}$ into $\text{Full}_{\mu,\lambda}$, $P\mapsto M^P$, such that
\begin{enumerate}
    \item for all $P\in \text{Full}_{\mu,\lambda}$, $M^P\subseteq P$,
    \item
    $M=M^{V_\lambda}$,
\item for all $P\in\text{Full}_{\mu,\lambda}$, $M^{M^P}=M^P$.
    \item  for any $P,Q\in  \text{Full}_{\mu,\lambda}$,
    if $\pi:H^P_{\mu^+}\to H^Q_{\mu^+}$ is elementary, then $\pi\restrict H^{M^P}_{\mu^+}$ is elementary from $H^{{M^P}}_{\mu^+}$ into $H^{{M^Q}}_{\mu^+}$,
   \item for any $P,Q\in  \text{Full}_{\mu,\lambda}$ such that $Q\subseteq P$, there is an elementary embedding $\pi:H^{M^P}_{\mu^+}\to H^{M^Q}_{\mu^+}$ such that
    \begin{enumerate}
        \item $\pi\in P$,
        \item $P\models``\pi$ is the unique elementary embedding from $H^{M^P}_{\mu^+}$ into $H^{M^Q}_{\mu^+}$".
    \end{enumerate}
\end{enumerate}
\end{definition}

Let us briefly discuss this definition. First, the function $P\mapsto M^P$ is really just proxy for $M$ being the output of some local definition of an inner model. This is why (4) is at all plausible. Still, it is convenient to abstract away from definability to the extent we can. Also note that (2) and (5) give that $H^M_{\mu^+}$ elementarily embeds into $H^{M^P}_{\mu^+}$ for any $P\in \text{Full}_{\mu,\lambda}$. Finally note that (5) for $P=Q=V_\lambda$ implies that there is no non-trivial elementary embedding $\pi:H^M_{\mu^+}\to H^M_{\mu^+}$, since the identity must be the unique such embedding (all such embeddings are in $V_\lambda$).

Under the hypothesis that there is no inner model with a Woodin cardinal, we will show that levels of $K$ resemble the core model, as witnessed by the function $P\mapsto (\varphi_K)^P$, and that the maps $\pi$ witnessing (5) for levels of $K$ are actually uniformly definable, which will be important for the uniqueness proof. We make the following definition capturing these additional properties we will verify for $K$.

\begin{definition}
Let $\mu$ be a measurable cardinal and $\lambda>\mu$ inaccessible. A transitive model $M$ \textit{strongly resembles the core model at $(\mu,\lambda)$} if there is a function $P\mapsto M^P$ such that (1)-(5) hold, $H^{M^P}_{\mu^+}$ is uniformly definable over $H^P_{\mu^+}$, and the maps $\pi:H^{M^P}_{\mu^+}\to H^{M^Q}_{\mu^+}$ witnessing (5) are definable over $H^P_{\mu^+}$, uniformly in parameter $H^Q_\mu$.
\end{definition}



We now prove the uniqueness result.

\begin{theorem}\label{uniqueness with measurables}
Suppose that $\mu$ is a measurable cardinal and $\lambda>\mu$ is inaccessible. Suppose that $N$ resembles the core model at $(\mu,\lambda)$ and $M$ strongly resembles the core model at $(\mu,\lambda)$. Then $H^M_{\mu^+}=H^N_{\mu^+}$.
\end{theorem}
\begin{proof}
Fix a function $P\mapsto N^P$ witnessing that $N$ resembles the core model at $(\mu, \lambda)$ and a function $P\mapsto M^P$ witnessing that $M$ \textit{strongly} resembles the core model at $(\mu,\lambda)$. We also fix $\varphi$ a formula witnessing that $H_{\mu^+}^{M^P}$ is uniformly definable over $H_{\mu^+}^P$ for $P\in \text{Full}_{\mu, \lambda}$, i.e. such that $H_{\mu^+}^{M^P}=\varphi^{H_{\mu^+}^P}$ for all $P\in \text{Full}_{\mu, \lambda}$ (such a $\varphi$ is guaranteed by the definition of strongly resembling the core model).

First we'll verify the following.
\begin{claim}
$H^{N^M}_{\mu^+} = H^M_{\mu^+}$\end{claim}

\begin{proof}
Let $M_0=M$, $N_0=N^{M}$, $M_1=M^{N_0}$, $N_1=N^{M_1}$, and $M_2=M^{N_1}$. Then, using (1) for $N$, $M_2\subseteq M_1\subseteq M_0$, so we can fix elementary embeddings $\pi:H^{M_0}_{\mu^+}\to H^{ M_1}_{\mu^+}$ and $\sigma:H^{M_1}_{\mu^+}\to H^{M_2}_{\mu^+}$ witnessing (5) for $M$. Similarly, fix an elementary embedding $\tau:H^{N_0}_{\mu^+}\to H^{N_1}_{\mu^+}$ witnessing (5) for $M$. (4) for $N$ gives $\pi\restrict H^{N_0}_{\mu^+}: H^{N_0}_{\mu^+}\to H^{N_1}_{\mu^+}$, so since $\pi\in M_0$ by (5)(a) for $M$, (5)(b) for $N$ gives that $\pi\restrict H^{N_0}_{\mu^+}=\tau$. A symmetric argument gives that $\tau\restrict H^{M_1}_{\mu^+}=\sigma$. So we have that $\pi\restrict H^{M_1}_{\mu^+}=\sigma$.

Now suppose that $\pi$ is not the identity and let $\kappa=\crit(\pi)$. Then $\kappa$ is definable over $H^{M_0}_{\mu^+}$ in parameter $H^{N_0}_\mu$, since $\pi$ is, using (3) for $M$ together with our assumption that $M$ \textit{strongly} resembles the core model at $(\mu,\lambda)$. These assumptions together with (4) for $N$ give that $\crit(\sigma)$ is defined in the same way over $H^{M_1}_{\mu^+}$ in parameter $H^{N_1}_\mu=\pi(H^{N_0}_\mu)$ as $\kappa=\crit(\pi)$ is over $H^{M_0}_{\mu^+}$ in parameter $H^{N_0}_\mu$. Since $\pi$ is elementary, it follows that $\pi(\kappa)=\crit(\sigma)$.  But $\sigma$ and $\pi$ agree on the ordinals, so $\crit(\sigma)=\crit(\pi)=\kappa$. So $\pi(\kappa)=\kappa$, contradicting that $\kappa$ is the critical point of $\pi$. So $\pi$ is the identity and $H^{N^M}_{\mu^+}=H^M_{\mu^+}$, as claimed \qed
\end{proof}

Next we show
\begin{claim}
$H^{M^N}_{\mu^+}=H^N_{\mu^+}$
\end{claim}
\begin{proof}
Since our hypotheses on $M$ and $N$ are not symmetric, this doesn't follow immediately from the proof of the previous claim. What that proof does give that $H^{N^{M^N}}_{\mu^+}=H^{M^N}_{\mu^+}$. So, by (5) for $N$, we get an elementary embedding $\pi:H^N_{\mu^+}\to H^{M^N}_{\mu^+}$. By (3) for $M$, we have that $H_{\mu^+}^{M^{M^N}}=H^{M^N}_{\mu^+}$.  In particular, $H^{M^N}_{\mu^+}\models \forall v\, \varphi(v)$, by our choice of $\varphi$. Since $\pi$ is elementary, $H^N_{\mu^+}\models \forall v\, \varphi(v)$ as well. So $H^{M^N}_{\mu^+}=H^{N}_{\mu^+}$, as claimed. \qed 
\end{proof}

By these claims and (5) for $M$ and $N$, there are elementary embeddings $\pi:H^M_{\mu^+}\to H^{M^N}_{\mu^+}=H^N_{\mu^+}$ and $\sigma: H^N_{\mu^+}\to H^{N^M}_{\mu^+}=H_{\mu^+}^M$. So $\sigma\circ \pi:H^M_{\mu^+}\to H^M_{\mu^+}$ is elementary and so must be the identity, by (5) for $M$ (see discussion following the definition). It follows that $\pi$ and $\sigma$ are the identity as well and so $H^M_{\mu^+}=H^N_{\mu^+}$, as desired.
\qed
\end{proof}

We now show that levels of $K$ strongly resemble the core model under the hypothesis that there is no inner model with a Woodin cardinal.

To start, we observe that the inductive definition of $K$ gives the following.

\begin{proposition}
Assume there is no inner model with a Woodin cardinal. Let $\mu$ be an inaccessible cardinal such that $\mu^{+,K}=\mu^+$. Then $K|\mu^+$ is definable without parameters over $H_{\mu^+}$, uniformly in $\mu$. 
\end{proposition}

\begin{proof}
 This follows from Theorem \ref{definability of K 1} together with the fact that, under the hypotheses of the proposition, $K|\mu^+=S(K|\mu)$, the stack of countably iterable sound premice extending $K|\mu$ which project to $\mu$. This is sufficiently definable by Theorem \ref{absoluteness2}, since the 1-small premice which are levels of this stack are cofinal. \qed
\end{proof}

The following is an immediate corollary to the previous proposition, Theorem \ref{definability of K 1}, and Theorem \ref{k of k}. 

\begin{corollary}
Assume there is no inner model with a Woodin cardinal. Let $\alpha$ be a limit cardinal of $K$ or the $K$-successor of an inaccessible cardinal of $K$. Let $P$ be a transitive model of ZFC${}^-$. Then
\begin{enumerate}
    \item if $\pi:H_\alpha^K\to P$ is elementary (in the language of set theory), there is a unique premouse $\hat{P}$ with universe $P$ such that $\pi:K|\alpha\to \hat{P}$ is elementary (in the language of premice), 
    \item if $\pi:P\to H_\alpha^K$ is elementary (in the language of set theory), there is a unique premouse $\hat{P}$ of $P$ such that is $\pi:\hat{P}\to K|\alpha$ is elementary (in the language of premice).
\end{enumerate}
\end{corollary}

Finally, we'll work towards showing that, in certain situations, elementary embeddings from initial segments of $K$ are uniquely determined by their target model. These results don't seem to appear in the literature, but are simple consequences of known techniques. 

We'll use the following easy criterion for being a fixed point of an embedding $\pi:M\to N$.
\begin{lemma}\label{fixed point lemma}
Let $M$, $N$ be transitive models of ZFC${}^-$ and $\pi:M\to N$ elementary. Suppose that $\sup\pi"\alpha=\alpha$ and $\pi$ is continuous at $\text{cof}^M(\alpha)$. Then $\pi(\alpha)=\alpha$.
\end{lemma}
\begin{proof}
Let $\gamma=\text{cof}^M(\alpha)$. Fix $\langle \beta_\xi\mid \xi<\gamma\rangle\in M$ a cofinal increasing sequence in $\alpha$. Then $\pi"\{ \beta_\xi\mid \xi<\gamma\}$ is cofinal in $\pi(\{ \beta_\xi\mid \xi<\gamma\})$ since $\sup\pi"\gamma=\pi(\gamma)$. So we have
\begin{align*}
    \pi(\alpha)&=\sup\pi(\{\beta_\xi\mid\xi<\gamma\})\\
   & = \sup\pi"\{\beta_\xi\mid\xi<\gamma\}\\
    &=\alpha.
\end{align*}\qed
\end{proof}

\begin{lemma}\label{cofinality lemma}
Let $\mu$ be a regular cardinal and $M$ transitive models of ZFC${}^-$ such that $o(M)=\mu^+$ and $\mu$ is the largest cardinal of $M$. Then for any $\alpha<\mu^+$, \[\text{cof}^M(\alpha)=\mu\Longleftrightarrow \text{cof}(\alpha)=\mu.\]
\end{lemma}
\begin{proof}
If an ordinal $\alpha<\mu^+$ has $\text{cof}(\alpha)=\mu$, then $\text{cof}^M(\alpha)=\mu$ since $\text{cof}(\mu)\leq\text{cof}^M(\alpha)\leq\mu$, as $|\alpha|^M\leq\mu$, since $\mu$ is the largest cardinal of $M$. Conversely, if $\text{cof}^M(\alpha)=\mu$, then $\text{cof}(\alpha)=\mu$ since $\mu$ is a regular cardinal.\qed
\end{proof}
\begin{definition}
Let $\mu$ be regular cardinal. We let $\mathcal{C}_{\mu,\mu^+}$ be the $\mu$-club filter on $\mu^+$; that is, the filter generated by the cofinal subsets of $\mu^+$ which are closed under increasing $\mu$-sequences.
\end{definition}

\begin{proposition}\label{fixed point 1.5}
Let $\mu$ be a regular cardinal and $M$, $N$ be transitive models of ZFC${}^-$ such that $o(M)=o(N)=\mu^+$ and $\mu$ is the largest cardinal of $M$. Suppose that $\pi:M\to N$ is elementary and $\pi$ is continuous at $\mu$. Then the set of fixed points of $\pi$ is a member of $\mathcal{C}_{\mu,\mu^+}$.
\end{proposition}
\begin{proof}
If $\alpha$ is a limit of fixed points of $\pi$ which has cofinality $\mu$, then Lemma \ref{cofinality lemma} gives $\text{cof}^M(\alpha)=\mu$, so that $\pi(\alpha)=\alpha$ by Lemma \ref{fixed point lemma} (since $\sup\pi"\alpha=\alpha$, as it is a limit of fixed points). So we just need to see $\pi$ has arbitrarily large fixed points. 

Fix $\beta<\mu^+$. Above $\beta$, we can build a $\mu$-sequence $\langle \alpha_\xi\mid \xi<\mu\rangle$  in $\mu^+$ such for all $\eta<\xi<\mu$, that $\pi(\alpha_\eta)<\alpha_{\xi}$. Let $\alpha=\sup\{\alpha_\xi\mid \xi<\mu\}$. Then $\alpha=\sup\pi"\alpha$ and, by Lemma \ref{cofinality lemma}, $\text{cof}^M(\alpha)=\text{cof}(\alpha)=\mu$. So by Lemma \ref{fixed point lemma}, $\pi(\alpha)=\alpha$. \qed
\end{proof}
We'll typically use this in the following situation. 
\begin{corollary}\label{fixed point 2}
Let $\mu$ be a regular cardinal and $M$, $N$ be transitive models of ZFC${}^-$ such that $o(M)=o(N)=\mu^+$ and $\mu$ is the largest cardinal of $M, N$. Suppose that $\pi:M\to N$ is elementary. Then the set of fixed points of $\pi$ is a member of $\mathcal{C}_{\mu,\mu^+}$.
\end{corollary}
\begin{proof}
Since $\mu$ is definable in the same way in $M,N$ (as the largest cardinal), $\pi(\mu)=\mu$. In particular, $\pi$ is continuous at $\mu$. So the previous proposition applies. \qed
\end{proof}

\begin{definition}
Let $\mu$ be a regular cardinal. An iterable premouse $P$ is \textit{$\mu$-universal} if $o(P)=\mu^+$ and $P$ has largest cardinal $\mu$.
\end{definition}

Note that, in general, there may be no premouse $P$ which is $\mu$-universal, according to this definition. However, if there is no inner model with a Woodin cardinal and there is a regular cardinal $\mu$ such that $\mu^{+,K}=\mu^+$ (e.g. for a  measurable cardinal $\mu$), then $K|\mu^+$ is $\mu$-universal.

\begin{theorem}\label{k at a measurable}
Assume there is no inner model with a Woodin cardinal. Suppose $\mu$ is a regular cardinal such that $\mu^{+,K}=\mu^+$. Let $P$ be $\mu$-universal. Then there is a unique elementary embedding $\pi:K|\mu^+\to P$. Moreover, $\pi$ is definable over $H_{\mu^+}$ in parameters $K|\mu$ and $P|\mu$, uniformly in $P|\mu$.
\end{theorem}

First we need to see that there is an embedding $\pi:K|\mu^+\to P$ at all. For this, we extend the definition of $\tilde{K}(\tau,\Omega)$ from \cite{kwm} to the case $\tau=\mu$ and $\Omega=\mu^+$.

\begin{definition}
Suppose that $P$ $\mu$-universal. $\text{Def}^P=\bigcap\{Hull^P(\Gamma)\mid \Gamma\in \mathcal{C}_{\mu, \mu^+}\}$.
\end{definition}

Standard arguments, as in \cite{kwm}, give

\begin{proposition}
Suppose that $P$ and $Q$ are $\mu$-universal. Then $\text{Def}^P\cong \text{Def}^Q$.
\end{proposition}

\begin{definition}
If there is a $\mu$-universal $P$, then $\tilde{K}(\mu, \mu^+)$ is the common transitive collapse of $\text{Def}^P$ for any $\mu$-universal $P$.
\end{definition} 

Now, the collapsing weasel case of the proof of Lemma 4.31 from \cite{kwm} gives
\begin{proposition}\label{k is local k prop}
Suppose that there is a $\mu$-universal $P$. Then $\tilde{K}(\mu, \mu^+)$ is $\mu$-universal and there is a $\Gamma\in \mathcal{C}_{\mu, \mu^+}$ such that $\text{Def}^P=\text{Hull}^P(\Gamma)$.
\end{proposition}

\begin{theorem}\label{k is local k}
Assume there is no inner model with a Woodin cardinal. Suppose $\mu$ is a regular cardinal such that $\mu^{+,K}=\mu^+$. Then $K|\mu^+=\tilde{K}(\mu, \mu^+)$.
\end{theorem}
\begin{proof}
We have that $K|\mu=\tilde{K}(\mu, \mu^+)|\mu$, since $\tilde{K}(\mu, \mu^+)|\mu$ also satisfies the inductive definition of $K$, as in the proof of Lemma 6.1 of \cite{kwm}. Since $\mu$ is regular, the stack over $K|\mu=\tilde{K}(\mu, \mu^+)|\mu$ is well-defined, i.e. the sound, iterable premice extending $K|\mu=\tilde{K}(\mu, \mu^+)|\mu$ and projecting to $\mu$ are totally ordered by the initial segment relation. It follows that $K|\mu^+\is \tilde{K}(\mu, \mu^+)$ or $\tilde{K}(\mu, \mu^+)\is K|\mu^+$. But both have height $\mu^+$, so they must be equal.
\qed
\end{proof}

This easily gives the following.

\begin{proposition}\label{k is local k prop2}
Assume there is no inner model with a Woodin cardinal. Suppose $\mu$ is a regular cardinal such that $\mu^{+,K}=\mu^+$. Then for any $\Gamma\in \mathcal{C}_{\mu,\mu^+}$, $K|\mu^+=\text{Hull}^{K|\mu^+}(\Gamma)$.
\end{proposition}
\begin{proof}
It is enough to show that $K|\mu^+=\text{Hull}^P(\Gamma)$ for \textit{some} $\Gamma\in \mathcal{C}_{\mu,\mu^+}$. Let $P$ be $\mu$-universal. Using Proposition \ref{k is local k prop} and Theorem \ref{k is local k}, we let $\Gamma$ be such that $K|\mu^+$ is the transitive collapse of $\text{Hull}^{P}(\Gamma)$. Let $\pi:K|\mu^+\to P$ be the uncollapse map. By Proposition \ref{fixed point 2}, we can assume that $\Gamma$ is a set of fixed points of $\pi$. For any $\Lambda\subseteq \Gamma$, since $\pi"\Lambda=\Lambda$, $\pi"\text{Hull}^{K|\mu^+}(\Lambda)=\text{Hull}^P(\Lambda)=\text{Def}^P$, and so $\ran(\pi)\subseteq \pi"\text{Hull}^{K|\mu^+}(\Lambda)$. It follows that $K|\mu^+=\text{Hull}^{K|\mu^+}(\Lambda)$. \qed
\end{proof}

\begin{proof}[Proof of Theorem \ref{k at a measurable}.]
By Theorem \ref{k is local k}, there is an embedding from $K|\mu^+$ into $P$, and by Proposition \ref{k is local k prop}, we actually have that $K|\mu^+$ is the transitive collapse of $\text{Def}^P=\hull^P(\Gamma)$ for some $\Gamma\in \mathcal{C}_{\mu,\mu^+}$. So suppose $\pi:K|\mu^+\to P$ is elementary. Then by Proposition \ref{fixed point 2}, the set of fixed points of $\pi$ is in $\mathcal{C}_{\mu,\mu^+}$, so we can find some $\Lambda\in \mathcal{C}_{\mu,\mu^+}$ which is a set of fixed points of $\pi$ such that $\text{Def}^P=\hull^P(\Lambda)$. We also get $K|\mu^+=\hull^{K|\mu^+}(\Lambda)$, by Proposition \ref{k is local k prop2}. It follows that
\begin{align*}
    \pi"K|\mu^+&=\pi"\hull^{K|\mu^+}(\Lambda)\\
    &=\hull^P(\pi"\Lambda)\\
     &=\hull^P(\Lambda)\\
     &=\text{Def}^P.
\end{align*} 
Since $\pi$ was arbitrary, $\text{Def}^P$ is the range of any elementary embedding from $K|\mu^+$ into $P$. So there is at most one such embedding. 

For the definability of $\pi$, since $\mu$ is regular, standard arguments give that $K|\mu$ and $P|\mu$ have a common iterate, $Q$, and letting $i:K|\mu\to Q$ and $j:P|\mu\to Q$ be the iteration maps of the comparison, and $E, F$ the length $\mu$ extenders of these iteration maps, \[S(Q)=\ult(K|\mu^+, E)=\ult(P, F),\] where $S(Q)$ is the stack over $Q$.
Let $\hat\imath:K|\mu^+\to S(Q)$ and $\hat\jmath:P\to S(Q)$ be the ultrapower maps. Then we also have that $\hat\imath"K|\mu^+=\hat\jmath"\text{Def}^P$. It follows that $\pi=\hat\jmath^{-1}\circ\hat\imath$. Since $K|\mu$, $P|\mu$, $Q$, and $E,F$ are all in $H_{\mu^+}$ and since $K|\mu^+=S(K|\mu)$ and $P=S(P|\mu)$, we get the required definability of $\pi$ (using for uniformity that $E,F$ came from the comparison). \qed
\end{proof}

\begin{lemma}
Assume there is no inner model with a Woodin cardinal. Suppose that $\mu$ is measurable and $\lambda>\mu$ is inaccessible. Then the function on $\text{Full}_{\mu,\lambda}$ given by $P\mapsto K^P=(\varphi_K)^P$ witnesses that $V_{\lambda}^K$ strongly resembles the core model.
\end{lemma}
\begin{proof}
First, $K^P$ is $\mu$-full at $\lambda$ since it is provable in ZFC $+$ ``there is no inner model with a Woodin cardinal" that for any measurable cardinal $\mu$, $\mu^{+,K}=\mu^+$ and $\mu$ is measurable in $K$. (This is part of Theorem \ref{k is close}, actually due to Steel alone.) So (1) holds. (2) follows trivially by how we chose our function. (3) follows from Theorem \ref{k of k}. (4), (5), and the additional definability requirement on the witnessing maps follow from Theorem \ref{k at a measurable} together with Theorem \ref{definability of K 1}. Finally, the fact that $H_{\mu^+}^{K^P}$ is definable over $H_{\mu^+}^P$ follows from the fact that, working inside $P$, $K|\mu$ is definable over $H_\mu$, since $\mu$ is a strong limit, by Theorem \ref{definability of K 1}, together with the fact that $K|\mu^+=S(K|\mu)$, the stack over $K|\mu$, as $\mu$ is regular and $\mu^{+,K}=\mu^+$. \qed
\end{proof}

This lemma and the previous theorem immediately imply the following.

\begin{theorem}
Assume there is no inner model with a Woodin cardinal. Suppose there is a proper class of measurable cardinals. Then $K$ is the unique inner model such that for all measurable cardinals $\mu<\lambda$, $V^K_\lambda$ resembles the core model at $(\mu,\lambda)$.
\end{theorem}

\section{Below $0^\text{\textparagraph}$}

In this section, we will prove a similar uniqueness result about the core model under the hypothesis that $0^\text{\textparagraph}$ does not exist. Here, we can drop our assumption that there is a proper class of measurable cardinals. Our approach will be the same as in the previous section: we will (re-)define what it means for an inner model to ``resemble the core model", prove that there is at most one such inner model, and prove that $K$ is this model. 

In this setting, the properties we'll use about the core model are known not just for the ms-indexed core model $K$ but for Schindler's core model below $0^\text{\textparagraph}$, which is defined via a different indexing scheme due to Jensen which we will call Jensen-indexing.\footnote{The indexing scheme used by Schindler is still different from the general indexing scheme for short extenders due to Jensen which is also called Jensen-indexing or sometimes $\lambda$-indexing. We thank Ralf Schindler for pointing this out to us.}  So our result gives that these two versions of the core model are the same under this hypothesis. Going forward, we denote by $J$ the Schindler core model below $0^\text{\textparagraph}$.

First, we'll need to look at directed systems of elementary embeddings between transitive models of ZFC, which we'll just call ``directed systems of models of set theory".

\begin{definition}
A \textit{directed system of models of set theory} is a system $\mathcal{D}=\{M_i, \pi_{i,j}\mid i,j\in D \text{ and }i\leq j\}$ such that \begin{enumerate}
    \item $\leq$ is a directed partial order on $D$,
    \item for every $i,j,k\in D$, \begin{enumerate}
    \item $M_i$ is a transitive model of ZFC,
    \item if $i\leq j$, then $\pi_{i,j}$ is an elementary embedding from $M_i$ into $M_j$,
    \item $\pi_{i,i}$ is the identity on $M_i$, and
    \item if $i\leq j\leq k$, then $\pi_{i,k}=\pi_{j,k}\circ \pi_{i,j}$.
    \end{enumerate}
\end{enumerate}
\end{definition}

\begin{definition}
For $\mathcal{D}=\{M_i, \pi_{i,j}\mid i,j\in D \text{ and }i\leq j\}$ a directed system of models of set theory,
$\mathcal{D}$ is \textit{well-founded} if the direct limit $(M_\infty, E)$ is well-founded, in which case we take $M_\infty$ to be transitive and $E$ to be $\in\restrict M_\infty$.
\end{definition}

We will use this terminology and the results to follow even when the directed system $\mathcal{D}$ is a definable family of transitive proper class models of ZFC (and elementary embeddings between them). Of course, as ZFC theorems, any results proven about such a system are schematic.

Our first lemma is implicit in the computations of HOD in models of determinacy, isolated in this general form by Gabriel Goldberg.

\begin{lemma}\label{elementlemma}
Let $\mathcal{D}=\{M_i, \pi_{i,j}\mid i,j\in D \text{ and }i\leq j\}$ be a well-founded directed system of models of set theory. Let $M_\infty$ be its direct limit, $\pi_{i,\infty}:M_i\to M_\infty$ the direct limit maps, and $X\subseteq M_\infty$.

The following are equivalent.
\begin{enumerate}
\item $X\in M_\infty$,
\item there is an $i\in D$ such that
\begin{enumerate}
    \item $\pi_{i,\infty}^{-1}[X]\in M_i$ and
    \item for any $j\geq i$, $\pi_{i,j}(\pi_{i,\infty}^{-1}[X])=\pi_{j,\infty}^{-1}[X]$.
\end{enumerate}
\end{enumerate}
\end{lemma}

\begin{proof}
First we show (1) $\Rightarrow$ (2). Suppose $X\in M_\infty$. Then $X$ is the image of an element of some model in our system, i.e. we can find an $i\in D$ and $\bar{X}\in M_i$ such that $\pi_{i,\infty}(\bar{X})=X$. We check $i$ witnesses (2) holds. For (2)(a), it's enough to see that $\bar{X}=\pi_{i,\infty}^{-1}[X]$. But this is trivial by elementarity: for any $x\in M_i$, 
\begin{align*}
    x\in \bar{X} &\Longleftrightarrow \pi_{i,\infty}(x)\in \pi_{i,\infty}(\bar{X})=X\\
    &\Longleftrightarrow x\in \pi_{i,\infty}^{-1}[X].
\end{align*}
Since for any $j\geq i$, $\pi_{j,\infty}(\pi_{i,j}(\bar X))=X$, the corresponding calculation at $j$ also gives (2)(b).

Now we show (2) $\Rightarrow$ (1). So let $i$ witness that (2) holds. For $j\geq i$, let $\bar{X}_j=\pi_{j,\infty}^{-1}[X]$. (2)(a) says $\bar{X}_i\in M_i$. Since $X\subseteq M_\infty$, (2)(b) gives that $\pi_{i,\infty}(\bar{X}_i)\subseteq M_\infty$, too. So, it's enough to show that $\pi_{i,\infty}(\bar{X}_i)$ and $X$ have the same elements of $M_\infty$. So fix $x\in M_\infty$.
Then $x$ is the image of an element of some point in our system, so we can find a $j\geq i$ and $\bar{x}\in M_j$ such that $\pi_{j,\infty}(\bar{x})=x$. By (2)(b), $\pi_{i,j}(\bar{X}_i)=\bar{X}_j$ (in particular, $\bar{X}_j\in M_j$). So since $\bar{X}_j=\pi_{j,\infty}^{-1}[X]$, 
\[\bar{x}\in \bar{X}_j \Longleftrightarrow x\in X.\]
Since $\pi_{j,\infty}(\bar{X}_j)=\pi_{i,\infty}(\bar{X}_i)$, applying $\pi_{j,\infty}$ to the left-hand side, gives
\[x\in \pi_{i,\infty}(\bar{X}_i) \Longleftrightarrow x\in X.\] \qed
\end{proof}

We'll identify a definability criterion which is sufficient for (2) and typical in applications of the lemma.

\begin{definition}
Let $\mathcal{D}=\{M_i, \pi_{i,j}\mid i,j\in D \text{ and }i\leq j\}$ be a directed system of models of set theory. For $A\subseteq D$, an $A$-indexed family of $n$-ary relations $\{R_i\mid i\in A\}$ is \textit{uniformly definable over $\mathcal{D}$} if there is an $i\in A$, $a\in M_i$, and formula in the language of set theory $\varphi(v_1,\ldots, v_n, u)$ such that for all $j\geq i$,\begin{enumerate}\item $j\in A$,
    \item
    $R_j\subseteq M_j^n$ and
    \item for all $x_1,\ldots, x_n\in M_j$, 
   $R_j(x_1,\ldots, x_n) \Longleftrightarrow M_j\models \varphi(x_1,\ldots, x_n, \pi_{i,j}(a))$.
\end{enumerate}
A single relation $R$ is \textit{uniformly definable over $\mathcal{D}$} if the constantly $R$ $D$-indexed family is uniformly definable over $\mathcal{D}$.
\end{definition}

\begin{definition}
For $M$ a transitive set, a set $X$ is a \textit{bounded subset of $M$} if there is a $y\in M$ such that $X\subseteq y$.
\end{definition}

\begin{lemma}\label{definabilitylemma1}
Let $\mathcal{D}=\{M_i, \pi_{i,j}\mid i,j\in D \text{ and }i\leq j\}$ be a well-founded directed system of models of set theory. Let $M_\infty$ be its direct limit, $\pi_{i,\infty}:M_i\to M_\infty$ the direct limit maps. Let $X$ be a bounded subset of $M_\infty$.

Suppose that $ \{ \pi_{i,\infty}^{-1}[X]\mid i\in D\}$ is uniformly definable over $\mathcal{D}$. Then $X\in M_\infty$.
\end{lemma}
\begin{proof}
Since $X$ is a bounded subset of $M_\infty$, we can fix $y\in M_\infty$ such that $X\subseteq y$. Let $i\in D$ and $\bar{y}\in M_i$ such that $\pi_{i,\infty}(\bar{y})=y$. Then for all $j\geq i$, $\pi_{j,\infty}[X]\subseteq \pi_{i,j}(\bar{y})$.  Since $ \{ \pi_{i,\infty}^{-1}[X]\mid i\in D\}$ is uniformly definable over $\mathcal{D}$, by increasing $i$ if necessary, we have that $\pi_{j,\infty}[X]$ is a bounded subset of $M_j$ which is definable over $M_j$, and so $\pi_{j,\infty}[X]\in M_j$ by Replacement in $M_j$. But then the elementarity of $\pi_{i,j}$ and the uniform definability of the $\pi_{i,\infty}^{-1}[X]$ immediately gives $\pi_{i,j}(\pi^{-1}_{i,\infty}[X])=\pi^{-1}_{j,\infty}[X]$, so (2)(b) holds as well. So $X\in M_\infty$ by Lemma \ref{elementlemma}.\qed
\end{proof}

We can use this lemma to show that appropriately intertwined directed systems of models of set theory have the same direct limit when the points in the systems are models of ``V=HOD" and the direct limit models and maps are uniformly definable over both systems. We state the result as a condition for when an initial segment of the direct limit of one system is a subset of another.

\begin{theorem}\label{directlimit thm1}
Let $\mathcal{C}$ and $\mathcal{D}$ be well-founded directed systems of models of set theory with the same underlying partial order, $\mathcal{C}=\{N_i, \sigma_{i,j}\mid i,j\in D \text{ and }i\leq j\}$ and $\mathcal{D}=\{M_i, \pi_{i,j}\mid i,j\in D \text{ and }i\leq j\}$. Let $N_\infty$, $M_\infty$ be the direct limit models and $\sigma_{i,\infty}$, $\pi_{i,\infty}$ be the direct limit maps.

Suppose that
\begin{enumerate}
        \item  for all $i\in D$, $N_i\subseteq M_i$,
    \item  for all $i\in D$, $N_i\models$ ``V=HOD", and
  
    \item  $N_\infty$, $M_\infty$, $\{\pi_{i,\infty}\}_{i\in D}$ are uniformly definable over $\mathcal{D}$.
\end{enumerate}
Then $N_\infty\subseteq M_\infty$.
\end{theorem}

\begin{proof}
Suppose that $N_\infty\not\subseteq M_\infty$. Since $N_i$ satisfies ``V=HOD", $N_\infty$ satisfies ``V=HOD", too. So we can look at the least set of ordinals $A\in N_\infty$ such that $A\not\in M_\infty$, under the definable well-order of $N_\infty$. Since $N_\infty$ and $M_\infty$ are uniformly definable over $\mathcal{D}$, so is $A$. 
Further, since $A$ is a member of $M_i$ for all sufficiently large $i$, by the uniform definability of the $N_\infty$, so $A$ is a bounded subset of $M_\infty$. Since $\{\pi_{i,\infty}\}_{i\in D}$ and $A$ are uniformly definable over $\mathcal{D}$, we easily get $\{\pi^{-1}_{i,\infty}[A]\mid i\in D\}$ is also uniformly definable over $\mathcal{D}$. So Lemma \ref{definabilitylemma1} gives $A\in M_\infty$, a contradiction. \qed
\end{proof}



Now we will review some results which hold below $0^{\text{\textparagraph}}$. Many of the relevant properties we've already mentioned for the ms-indexed core model $K$ are known to hold for Schindler's core model $J$. The basic theory of this model is developed in Schindler's \cite{jensen-indexed core model} (see also \cite{zeman}).

The covering theorems of Schindler \cite{jensen-indexed core model} and Cox \cite{cox} immediately give
\begin{theorem}\label{j is close}
Assume $0^{\text{\textparagraph}}$ does not exist. Then $J$ is close to $V$.
\end{theorem}

Schindler's core model $J$ also has an inductive definition which gives a version of Theorem \ref{definability of K 1} under the hypothesis that $0^{\text{\textparagraph}}$ does not exist. We let $\varphi_J$ be the formula defining $J$ in this way. We also let $``V=J"$ be the formula $\forall v\, \varphi_J(v)$.

One important feature of inner model theory below $0^{\text{\textparagraph}}$ is that the theory of (definable) proper class premice is well-behaved. This is because definable iteration trees on iterable premice have definable well-founded branches, even when they are proper class sized. This fails below a Woodin cardinal in general. 

We have the following below $0^{\text{\textparagraph}}$.

\begin{theorem}\label{definable maps thm}
Assume $0^{\text{\textparagraph}}$ does not exist. Suppose that $\varphi$ is a $\Sigma_n$-formula which defines a close inner model $W$. Then
\begin{enumerate}
    \item there is an elementary embedding $k:K\to K^W$, definable uniformly in $\varphi$,
    \item there is an elementary embedding $j:J\to J^W$, definable uniformly in $\varphi$.
\end{enumerate}
\end{theorem}

\begin{proof}[Proof sketch.]
We'll just talk about the ms-indexed core model $K$, as it is symmetric. $K^W$ is iterable in $W$ and so it is actually iterable, by Theorem \ref{absoluteness}.\footnote{Below $0^\text{\textparagraph}$, this absoluteness fact is actually easier and holds for Jensen-indexed premice as well.} Moreover, $K^W$ is close to $M$ by Theorem \ref{k is close} in $W$. So $K^M$ is close to $V$ by Proposition \ref{transitivity}. In particular, $(\mu^+)^{K^M}=\mu^+$ for all singular strong limit cardinals $\mu$, so $K^M$ is \textit{universal} in the sense that it is maximal in the ms-indexed mouse-order, by standard arguments (cf. Lemma 6.3.1 of Zeman \cite{zeman}). In particular, $K$ and $K^M$ have a common, non-dropping iterate, obtained by comparing the two inner models. By standard arguments (which can be found in \cite{zeman} or \cite{cmip}), $K^W$ doesn't move in this comparison and so there is an elementary embedding $k:K\to K^W$. Since $k$ was obtained as the iteration map of the (definable) comparison, it is definable, uniformly in the definition of $W$.
\qed
\end{proof}
We suspect that Theorem \ref{definable maps thm} may fail below a Woodin cardinal, but we do not have a counterexample.

Finally, we also have that $K$ and $J$ are \textit{rigid}. Because rigidity of an inner model is not expressible in the language of set theory, in general, we make the following definition.

\begin{definition}
An inner model $M$ is \textit{$\Sigma_n$-rigid} if there is no $\Sigma_n$-definable, non-trivial elementary embedding $j:M\to M$.
\end{definition}

By standard techniques (cf \cite{zeman} or \cite{cmip}), we have the following, for any $n$.

\begin{theorem}\label{rigidity thm}
Assume $0^{\text{\textparagraph}}$ does not exist. Then
\begin{enumerate}
    \item $K$ is $\Sigma_n$-rigid,
    \item $J$ is $\Sigma_n$-rigid.
\end{enumerate}
\end{theorem}

This easily implies the following analogue of Theorem \ref{k of k}.

\begin{theorem} \label{j of j}
Assume $0^{\text{\textparagraph}}$ does not exist. Then $J\models ``V=J"$. In particular, $J\models ``V=HOD"$.
\end{theorem}
\begin{proof}[Proof sketch.]
 By Theorem \ref{definable maps thm} in $V$, there is a definable elementary embedding $j_0:J\to J^J$. But also since $J$ is a universal Jensen-indexed proper class premouse (i.e. Jensen-indexed mouse-order maximal), the proof of \ref{definable maps thm}\footnote{We mean the proof of the Jensen-indexed case, which was omitted but is symmetric to the ms-indexed case, which we sketched.} gives there is also an elementary embedding $j_1:J^J\to J$ which is still definable in $V$. But then $j_1\circ j_0:J\to J$ and so must be the identity, by Theorem \ref{rigidity thm}. It follows that $J=J^J$, as desired. We get $J\models ``V=HOD$" just as in the proof of Theorem \ref{k of k}. \qed
\end{proof}

Of course, this same argument gives another proof of Theorem \ref{k of k} under the additional hypothesis that $0^{\text{\textparagraph}}$ does not exist.

In \cite{cmip}, Steel proves that if there is no inner model with a Woodin cardinal and $\Omega$ is a measurable cardinal, then $K|\Omega$ is rigid. Surprisingly, it appears to be open whether $K$ is ($\Sigma_n$-)rigid just under the hypothesis that there is no inner model with a Woodin cardinal. The difficulty in adapting arguments from \cite{cmip}  to this context is that these arguments rely on the existence of very soundness witnesses for initial segments of $K$ which are definable over $K$. We do not see how to get such witnesses in the general context.



Next, we'll look at inner models which are definable via special kinds of $\Sigma_2$-formulas, where by \textit{inner models} we mean transitive proper class models of ZF, as is standard. We will look at formulas which provably define inner models over some base theory, $T$. Our main interest is in the theory ZFC$+``0^{\text{\textparagraph}}$ does not exist". The hypothesis that $0^{\text{\textparagraph}}$ does not exist are both absolute to inner models: that is, if it holds, it holds in any inner model of ZFC. We introduce the following definition to capture this phenomenon.

\begin{definition}
A theory $T$ is a \textit{nice} extension ZFC iff it has the form ZFC$+\varphi$ for a $\Pi_2$-sentence $\varphi$ of the form ``for every set of ordinals $A$, $L[A]\models \theta$" for some $\Pi_2$-sentence $\theta$. \end{definition}

We leave the following easy proposition to the reader.
\begin{proposition} Let $T$ be a nice extension of ZFC. Assume $T$. Then for any transitive proper class model $W$ of ZFC, $W\models T$.
\end{proposition}

\noindent It is straightforward to see that ZFC$+``0^{\text{\textparagraph}}$ does not exist" is nice extensions of ZFC.

We'll use the following standard notation. \begin{definition}
For $\varphi(\vec{u},\vec{v})$ a formula with free variables $\vec{u},\vec{v}$, $M$ a transitive class and $\vec{x}\in M$, we let \[\varphi(\vec{x}, \vec{v})^M=\{\vec y\in M\mid M\models \varphi(\vec{x},\vec{y})\}.\]
Also, if $\varphi(u)$ has just one free variable $u$, we'll write $\varphi^M$ instead of $\varphi(u)^M$.
\end{definition}

\begin{definition}
Let $T$ be a nice extension of ZFC. A $\Sigma_2$-formula $\varphi(v)$ \textit{locally defines an inner model over $T$} iff $\varphi(v)$ has the form
\[\exists \mu\,\big(\mu \text{ is a strong limit cardinal } \wedge v\in H_{\mu}\wedge H_{\mu}\models \theta(v)\big).\]
for some formula $\theta(v)$ and, letting $M=\varphi^V$, the following is provable over $T$:
\begin{itemize}
    \item $M$ is an inner model of ZFC,\footnote{Recall that this is expressible in the language of set theory by asserting that $\varphi^V$ is a transitive class which is almost universal, closed under G\"odel operations, and satisfies AC; see \cite{Jech}. Note that since we assumed $T$ is a nice extension of ZFC, it follows that $M\models T$, as well.}
    \item for every strong limit cardinal $\mu$, $H^M_\mu=\theta^{H_{\mu}}$.
\end{itemize}
\end{definition}
\noindent Note that if $\varphi$ locally defines an inner model over $T$, then $\varphi$ is a $\Sigma_2$-formula. Also note that we can always take the $\mu$ in the displayed formula above to be the least strong limit cardinal such that $v\in H_{\mu}$.
We'll see that these formulas are more nicely behaved than arbitrary $\Sigma_2$-formulas which provably define inner models.

We need one more bit of notation.

\begin{definition}
Let $\rho=\rho(u_1,\ldots, u_n)$ and $\chi=\chi(\vec{v}, w)$ be formulas. The formula $\rho^\chi$ is \[{\rho^\chi}^-\wedge\chi(\vec{v}, u_1)\wedge\cdots \wedge \chi(\vec{v}, u_n),\] where ${\rho^\chi}^-$ defined recursively on the complexity of $\rho$ as follows:
\begin{itemize}
    \item for $\rho$ an atomic formula, ${\rho^\chi}^-= \rho$,
    \item ${(\rho\wedge\xi)^\chi}^-={\rho^\chi}^-\wedge {\xi^\chi}^-$,
    \item ${(\rho\vee\xi)^\chi}^-={\rho^\chi}^-\vee {\xi^\chi}^-$,
    \item ${(\neg\rho)^\chi}^-=\neg({\rho^\chi}^-)$,
    \item ${(\exists u\,\rho)^\chi}^-=\exists u \,(\chi(\vec{v},u)\wedge{\rho^\chi}^-)$, and
    \item ${(\forall u\,\rho)^\chi}^-=\forall u\, \big(\chi(\vec{v},u)\to {\rho^\chi}^-\big)$.
\end{itemize}
\end{definition}

The point of this is just that if $\chi^V$ is a transitive class, say $\chi^V=M$, then $(\rho^\chi)^V=\rho^M$.

\begin{lemma}\label{local1}
Let $T$ be a nice extension of ZFC. Let $\varphi,\psi$ be formulas which locally define inner models over $T$. Then there is a formula $\tau$ which locally defines an inner model over $T$ such that $\tau$ equivalent to $\psi^\varphi$, provably over $T$.
\end{lemma}
\begin{proof}
Let $\theta, \rho$ be formulas witnessing that $\varphi$ and $\psi$ locally define inner models over $T$, that is, such that $\varphi$ is $\exists \mu\,\big(\mu \text{ is a strong limit cardinal } \wedge v\in H_{\mu}\wedge H_{\mu}\models \theta(v)\big)$ and $\psi$ is
$\exists \mu\,\big(\mu \text{ is a strong limit cardinal } \wedge v\in H_{\mu}\wedge H_{\mu}\models \rho(v)\big)$.

Let $\tau$ be $\exists \mu\,\big(\mu \text{ is a strong limit cardinal } \wedge v\in H_{\mu}\wedge H_{\mu}\models \rho^\theta(v)\big)$.
We'll show that $\tau$ is our desired formula. Work in $T$. Let $M=\varphi^V$ and $N=\psi^M$. Then $M$ and $N$ are both inner models of $T$. We need to show that $N=\tau^V$ and that for any strong limit $\mu$, $H^N_{\mu}=(\rho^\theta)^{H_{\mu}}$. This latter claim immediately implies the former, so we just need to verify it. 

Let $\mu$ be a strong limit cardinal. Then $\mu$ is a strong limit cardinal of $M$, so \begin{align*}
    H_{\mu}^N&=\rho^{H^M_{\mu}}\\
    &=\rho^{(\theta^{H_{\mu}})}\\
    &=(\rho^\theta)^{H_{\mu}},
\end{align*}
using that $\varphi$ and $\psi$ locally define inner models over $T$ (as witnessed by $\theta$ and $\rho$) for the second and first equivalences, respectively. \qed
\end{proof}

For arbitrary $\Sigma_2$-formulas $\psi$ and $\varphi$ which provably define inner models, it seems that $\psi^\varphi$ should not be provably equivalent to a $\Sigma_2$-formula, but we do not have an example.


Our next goal is to shows that for $M, N$ inner models defined via local formulas over $T$, the $\omega$-sequence of inner models $\langle M, N^M, M^{N^M}, N^{M^{N^M}},\ldots\rangle$ is definable. The problem is that this we may have no bound on the quantifier complexity of the formulas $\varphi,\varphi^\psi, \psi^{\varphi^\psi}, \varphi^{\psi^{\varphi^\psi}},\ldots$ (where $\varphi$, $\psi$ are some witnessing formulas to the definability of $M,N$). We get around this by using our previous proposition.

Fix $\langle \varphi_i\mid i\in \omega\rangle$ a primitive recursive enumeration of formulas of the language of set theory in one free variable. For $\gamma$ such a formula, let $\ulcorner \gamma\urcorner$ be the $i$ such that $\gamma=\varphi_i$. For $\varphi,\psi$ formulas, let $F_{\varphi,\psi}$ be the primitive recursive function outputting the G\"odel numbers of the sequence $\langle \varphi,\psi^\varphi,\psi^{\varphi^\psi},\ldots\rangle$. That is, $F_{\varphi,\psi}$ is the function $F$ defined by
\begin{itemize}
    \item $F(0)=\ulcorner \varphi\urcorner$ and \item$F(k+1)=\begin{cases}\ulcorner\psi^{\varphi_{F(k)}}\urcorner & \text{if $k$ is even,}\\
    \ulcorner\varphi^{\varphi_{F(k)}}\urcorner & \text{if $k$ is odd.}\end{cases}$
\end{itemize}Let $\text{Sat}_0(w,v,u)$ be the usual definition of the $\Delta_0$-satisfaction predicate.\footnote{So for any transitive $x$, $y\in x$, and $n\in\omega$, $\text{Sat}_0(x,n,y)\Leftrightarrow x\models \varphi_n(y)$.} 

For $\varphi$, $\psi$ formulas which locally define inner models over $T$, say $\varphi$ is $\exists \mu\,\big(\mu \text{ is a strong limit cardinal } \wedge v\in H_{\mu}\wedge H_{\mu}\models \theta(v)\big)$ and $\psi$ is
$\exists \mu\,\big(\mu \text{ is a strong limit cardinal } \wedge v\in H_{\mu}\wedge H_{\mu}\models \rho(v)\big)$,
 we also let $\xi_{\varphi,\psi}(u,v)$ be the formula
\[\exists \mu \big( \mu \text{ is a strong limit cardinal }\wedge v\in H_{\mu}\wedge \text{Sat}_0(H_{\mu}, F_{\theta,\rho}(u), v)\big).\footnote{Here we really mean that we've replaced $F_{\theta,\rho}$ with a formula defining defining it over ZFC.}\]

\begin{proposition}\label{localmain2}
Let $T$ be a nice extension of ZFC. Suppose that $\varphi,\psi$ are formulas which locally define inner models over $T$. Let $\xi=\xi_{\varphi,\psi}$. Then the following is provable in $T$.\\

For every $k\in \omega$,
\begin{enumerate}
    \item $\xi(k,v)^V$ is an inner model of $T$\footnote{Here we mean $\xi(k,v)^V$ is almost universal, etc., and satisfies the additional sentence witnessing that $T$ is nice.}
        \item $\xi(0,v)^V=\varphi^V$,
        \item $\xi(k+1,v)^V=\psi^{\xi(k,v)^V}$ if $k$ is even, and
        \item $\xi(k+1,v)^V=\varphi^{\xi(k,v)^V}$ if $k$ is odd.
\end{enumerate}
\end{proposition}

\begin{proof}
Assume $T$. We prove (1)-(4) by induction on $k\in\omega$.

Since $F_{\theta,\rho}(0)=\ulcorner \theta\urcorner$, we immediately get $\xi(0,v)^V=\varphi^V$, giving (2). Since $\varphi^V$ is an inner model of $T$, by hypothesis, (1) holds for $k=0$.

Now suppose (1) holds at $k$, i.e. $\xi(k,v)^V$ is an inner model of $T$. Assume $k$ is even. We'll verify $\xi(k+1,v)^V=\psi^{\xi(k,v)^V}$. Then, since $\xi(k,v)^V$ is an inner model of ZFC, so is $\xi(k+1,v)^V=\psi^{\xi(k,v)^V}$ by our hypothesis about $\psi$. We have that $\xi(k, v)$ is
\[\exists \mu \big( \mu\text{ is a strong limit cardinal}\wedge v\in H_{\mu}\wedge \text{Sat}_0(H_{\mu}, F_{\theta,\rho}(k), v)\big),\]
which is equivalent to
\[\exists \mu \big( \mu\text{ is a strong limit cardinal}\wedge v\in H_{\mu}\wedge H_{\mu}\models\varphi_{F_{\theta,\rho}(k)}(v)\big).\]
Now, since $k$ is even,
$F_{\theta,\rho}(k+1)=\ulcorner \rho^{\varphi_{F_{\theta,\rho}}}\urcorner$. So, $\xi(k+1, v)$ is equivalent to
\[\exists \mu \big( \mu\text{ is a strong limit cardinal}\wedge v\in H_{\mu}\wedge H_{\mu}\models \rho^{\varphi_{F_{\theta,\rho}(k)}}( v)\big).\]
By the proof of Lemma \ref{local1}, we get that $\xi(k+1,v)^V=\psi^{\xi(k,v)^V}$, as desired. The case that $k$ is odd is basically the same (just replace $\rho$ with $\theta$).
\qed
\end{proof}

Again, although the precise statement is technical, we think of this proposition as saying that for $M, N$ locally definable inner models (over some nice $T$), the $\omega$-sequence of inner models $\langle M, N^M, M^{N^M}, N^{M^{N^M}},\ldots\rangle$ is actually definable (over $T$). This proposition is a major reason why we've focused on inner models which are locally definable: it is not clear that the sequence $\langle M, N^{M}, M^{N^{M}},\ldots\rangle $ is definable at all for two arbitrary $\Sigma_2$-formulas which provably define inner models (over some $T$), as the quantifier complexity of the resulting definitions gets arbitrarily large. 

We make another definition which is just a strengthening of locally defining an inner model over $T$.

\begin{definition}
A formula $\varphi$ \textit{locally defines a close inner model over $T$} iff $\varphi$ provably defines a close inner model over $T$ and $T$ proves that $\varphi^V$ is close.
\end{definition}

The most important examples of such formulas are $\varphi_K$ and $\varphi_J$ which both locally define close inner models over the theory ZFC$+``0^{\text{\textparagraph}}$ does not exist". This follows from our observations about $K$ in the prelimary section and the preceding observations about $J$.

Now we can state our main definition.
\begin{definition}
Let $T$ be a nice extension of ZFC. For $M$ an inner model, $\varphi(v)$ a formula which locally defines a close inner model over $T$, and $\psi(u,v,w)$ a formula, \textit{$M$ resembles the core model via }$(\varphi,\psi)$ iff $M=\varphi^V$ and the following is provable in $T$:
\begin{enumerate}
    \item $\varphi^V\models $``V=HOD",
    \item $\varphi^V\models\forall x\,\varphi(x)$,
    \item for any $\Sigma_2$-formula $\gamma$, if $W=\gamma^V$ is a close inner model of ZFC, then $\psi(\ulcorner\gamma\urcorner,v,w)$ defines an elementary embedding from $\varphi^V$ into $\varphi^W$;\footnote{Recall that an elementary embedding between definable inner models of ZFC is expressible in the language of set theory by asserting just $\Sigma_1$-elementarity.} 
\end{enumerate}
We'll say that $M$ $\Sigma_n$-\textit{resembles the core model over $T$} if we can take $\psi$ to be a $\Sigma_n$-formula. We'll just say that $M$ \textit{resembles the core model} if it $\Sigma_{10000}$-resembles the core model over ZFC$+``0^{\text{\textparagraph}}$ does not exist". Of course, this number is overkill: we just need $n$ sufficiently large so that the actual core model, $K$, $\Sigma_n$-resembles the core model over ZFC$+``0^{\text{\textparagraph}}$ does not exist".
\end{definition}

Our main result is the following schema, for $n\geq 1$.

\begin{theorem}\label{main 1}
Assume $T$. Suppose $M$ and $N$ resemble the core model over $T$ via $\Sigma_n$-formulas and $M$ is $\Sigma_n$-rigid. Then $M=N$.
\end{theorem}

\begin{proof}
Fix $\varphi,\rho$ such that $M$ resembles the core model via $(\varphi, \rho)$ and let $\psi, \chi$ such that $N$ resembles the core model via $(\psi, \chi)$. So $\varphi$ and $\psi$ are formulas which locally define close inner models over $T$. For $W$ an inner model, we'll write $M^W$ instead of $\varphi^W$ and $N^W$ instead of $\psi^W$.

We first show $N^M=M$ and $M^N=N$. Since it's symmetric, we'll just show the former. To get this, we'll show that $M$ elementarily embeds into some $M_\infty$ such that $M_\infty\models \forall x\, \psi(x) $, so that $M\models \forall x\, \psi(x) $, i.e. $N^M=M$.

Let $\xi=\xi_{\varphi,\psi}$ be the formula from Proposition \ref{localmain2} (defined in the discussion preceeding it). Let $M_k=\xi(2k, v)^V$ and $N_k=\xi(2k+1, v)^V$. Then by Proposition \ref{localmain2}, $M_0=\varphi^V=M$, $N_k=N^{M_k}$, and $M_{k+1}=M^{N_k}$. Let $\pi_i=\rho(\ulcorner \psi\urcorner,u,v)^{M_i}$ and $\sigma_i=\chi(\ulcorner \varphi\urcorner, u, v)^{N_i}$. So $\pi_i$ is an elementary embedding from $M_i$ into $M_{i+1}$ and $\sigma_i$ is an elementary embedding from $N_i$ into $N_{i+1}$. For $i\leq j$, let $\pi_{i,j}:M_i\to M_j$ and $\sigma_{i,j}:N_i\to N_j$ the natural maps obtained from composing the $\pi_k$ and $\sigma_k$, respectively. Let $\mathcal{C}=\{ N_i, \sigma_{i,j}\mid i,j\in\omega\text{ and } i\leq j\}$ and $\mathcal{D}=\{ M_i, \pi_{i,j}\mid i,j\in\omega\text{ and } i\leq j\}$. Let $N_\infty$ be the direct limit of $\mathcal{C}$, $M_\infty$ the direct limit of $\mathcal{D}$, and $\sigma_{i,\infty}$ and $\pi_{i,\infty}$ the direct limit maps.

First we'll show
$M_\infty, N_\infty$ are well-founded.
Of course, the argument is the same for $M_\infty$ and $N_\infty$. In fact, the argument is just Gaifman's argument showing that the $\omega^{\text{th}}$-iterate of $V$ by a countably complete ultrafilter is well-founded.

We'll just show $M_\infty$ is well-founded. To run Gaifman's argument, the main thing we need to observe is that $\{\langle M_{i+j}\mid j<\omega\rangle\rangle\mid i<\omega\}$ and $\{\langle\pi_{i+j,i+k}\mid j\leq k<\omega\rangle \mid i< \omega\}$ are both uniformly definable over $\mathcal{D}$. That $\{\langle M_{i+j}\mid j<\omega\rangle\rangle\mid i<\omega\}$ is uniformly definable over $\mathcal{D}$ is immediate since $M_{i+j}= \xi(2j, v)^{M_i}$ for all $i,j<\omega$, by Proposition \ref{localmain2}. That $\{\langle\pi_{i+j,i+k}\mid j\leq k<\omega\rangle \mid i< \omega\}$ is uniformly definable over $\mathcal{D}$ follows since $\pi_{i+j, i+k}$ is just the composition of $\rho(\ulcorner \psi\urcorner, u, v)^{M_{i+k-1}}\circ \cdots \circ\rho(\ulcorner \psi\urcorner, u, v)^{M_{i+j}}$. (That this is uniformly definable uses that the sequence of models $M_{i+j}$, $\ldots$, $M_{i+k-1}$ is uniformly definable, which follows from our preceding observation.)

Now suppose $M_\infty$ is not well-founded. Let $\alpha$ least such $M_\infty$ is ill-founded below $\pi_{0,\infty}(\alpha)$. But then by the uniform definability of  $\{\langle M_{i+j}\mid j<\omega\rangle\rangle\mid i<\omega\}$ and $\{\langle\pi_{i+j,i+k}\mid j\leq k<\omega\rangle \mid i< \omega\}$ are both uniformly definable over $\mathcal{D}$, for any $i$, $\pi_{0,i}(\alpha)$ is the least $\beta$ such that $M_\infty$ is ill-founded below $\pi_{i,\infty}(\beta)$. So $M_\infty$ cannot be ill-founded below $\pi_{0,\infty}(\alpha)$ after all, a contradiction.

As we mentioned in the preceding argument, $\{\langle\pi_{i+j,i+k}\mid j\leq k<\omega\rangle \mid i< \omega\}$ is uniformly definable over $\mathcal{D}$. So by Replacement, the (transitivized) $M_\infty$ and a tail of the direct limit maps $\pi_{i,\infty}$ are also uniformly definable over $\mathcal{D}$. Since $N_i$ is definable in $M_i$ we also get the (transitivized) $N_\infty$ and the direct limit maps $\sigma_{i,\infty}$ are also uniformly definable over $\mathcal{D}$. Since we also have $N_i\subseteq M_i$ and $N_i\models$ ``V=HOD" (by clause (1) of the definition of resembles the core model), Theorem \ref{directlimit thm1} gives $N_\infty\subseteq M_\infty$. A symmetric argument shows $M_\infty\subseteq N_\infty$. But by elementarity, $N_\infty\models \forall x\, \psi(x)$ so $M_\infty$ does too, as desired. 

To finish, let $\pi=\rho(\ulcorner \psi\urcorner,u,v)^V$ and $\sigma=\chi(\ulcorner \varphi\urcorner,u,v)^V$. So $\pi:M\to M^N=N$ and $\sigma:N\to N^M=M$ are elementary. If $M\neq N$ then at least one of $\pi,\sigma$ is not the identity on the ordinals. But then as $\sigma, \pi$ are definable by the $\Sigma_n$-formulas $\psi$ and $\tau$, $\sigma\circ \pi:M\to M$ is a $\Sigma_n$-definable elementary embedding\footnote{Using here that $n\geq 1$.} which is not the identity, contradicting the $\Sigma_n$-rigidity of $M$. So $M=N$ after all.
\qed
\end{proof}

Recall that $K$ is ms-indexed Jensen-Steel core model and $J$ is Schindler's core model below $0^\text{\textparagraph}$.

\begin{theorem}\label{maintheorem} Assume $0^\text{\textparagraph}$ doesn't exist. Then $K=J$ is the unique inner model which resembles the core model.
\end{theorem}

\begin{proof}
The inductive definitions of $K$ and $J$, given by formulas $\varphi_K$ and $\varphi_J$, locally define close inner models over ZFC$+``0^\text{\textparagraph}$ doesn't exist". We want to see that there are formulas $\psi_K$ and $\psi_J$ such that $K$ and $J$ resemble the core model via $(\varphi_K, \psi_K)$ and $(\varphi_J, \psi_J)$, respectively. But (1) and (2), which only mention $\varphi$, are immediate from Theorems \ref{k of k} and \ref{j of j}. Also, Theorem \ref{definable maps thm} immediately gives us our desired formulas $\psi_K$ and $\psi_J$ witnessing (3) for $\varphi_K$ and $\varphi_J$. Finally, we actually have that both $J$ and $K$ are sufficiently rigid, by Theorem \ref{rigidity thm}, so Theorem \ref{main 1} gives $K=J$ is the unique inner model which resembles the core model. \qed
\end{proof}


\begin{thebibliography}{}

\bibitem{cox} Sean Cox, Covering theorems for the core model, and an application to stationary set reflection, \textit{Annals of Pure and Applied Logic}, 161, 2009.

\bibitem{fuchs1} Gunter Fuchs, $\lambda$-structures and s-structures: Translating the models, \textit{Annals of Pure and Applied Logic}, 162 (4), 2011.

\bibitem{fuchs2} Gunter Fuchs, $\lambda$-structures and s-structures: Translating the iteration strategies, \textit{Annals of Pure and Applied Logic}, 162 (9), 2011.

\bibitem{schindler} Moti Gitik, Ralf Schindler, and Saharon Shelah, Pcf theory and Woodin cardinal, \textit{Logic Colloquium '02}, \textit{Lecture Notes in Logic}, 27, 2006.

\bibitem{Jech} Thomas Jech, \textit{Set Theory: the Third Millennium Edition}, Springer, 2003. 

\bibitem{jensen manu} Ronald B. Jensen, Manuscript on fine structure, inner model theory, and the core model below one Woodin cardinal, Draft, 2020.

Available at \url{https://www.mathematik.hu-berlin.de/~raesch/org/jensen.html}


\bibitem{kwm} Ronald B, Jensen and John R. Steel, K without the measurable, \textit{The Journal of Symbolic Logic}, 73, 2013.

\bibitem{mitchell-schimmerling} William J. Mitchell and Ernest Schimmerling, Covering at limit cardinals of K, \textit{Journal of Mathematical Logic}, 24 (01), 2024.

\bibitem{jensen-indexed core model} Ralf Schindler, The core model for almost linear iterations, \textit{Annals Pure and Applied Logic}, 116, 2002.

\bibitem{schindler steel}
Ralf Schindler and John R. Steel, The self-iterability of $L[E]$. \textit{Journal of Symbolic Logic}, 74, 2009.

\bibitem{thesis} Benjamin Siskind, \textit{Aspects of Martin's Conjecture and Inner Model Theory}. Ph.D. thesis,
UC Berkeley, 2021.
\bibitem{cmip} John R. Steel, \textit{The Core Model Iterability Problem}. Springer-Verlag, Berlin, 1996.


\bibitem{pfa} John R. Steel, PFA implies AD in $L(\mathbb{R})$. \textit{Journal Symbolic Logic}, 70 (4), 2005.


\bibitem{steel book}
John R. Steel, \textit{A Comparison Process for Mouse Pairs}, Lecture Notes in Logic, Cambridge University Press, 2024.

\bibitem{hod_as}
John R. Steel and W. Hugh Woodin, HOD as a core model, In: Kechris AS, Löwe B, Steel JR, eds. \textit{Ordinal Definability and Recursion Theory: The Cabal Seminar, Volume III}, Lecture Notes in Logic, Cambridge University Press, 2016.



\bibitem{zeman} Martin Zeman, \textit{Inner Models and Large Cardinals}. \textit{de Gruyter Series in Logic and Applications}, 5, de Gruyter, Berlin, 2002.




\end{thebibliography}
\end{document}